\documentstyle[12pt]{article}
\begin{document}
\renewcommand{\thefootnote}{\fnsymbol{footnote}}
\newpage
\pagestyle{empty}
\setcounter{page}{0}
\renewcommand{\thesection}{\arabic{section}}
\renewcommand{\theequation}{\thesection.\arabic{equation}}
\newcommand{\sect}[1]{\setcounter{equation}{0}\section{#1}}
\newfont{\twelvemsb}{msbm10 scaled\magstep1}
\newfont{\eightmsb}{msbm8}
\newfont{\sixmsb}{msbm6}
\newfam\msbfam
\textfont\msbfam=\twelvemsb
\scriptfont\msbfam=\eightmsb
\scriptscriptfont\msbfam=\sixmsb
\catcode`\@=11
\def\Bbb{\ifmmode\let\next\Bbb@\else
  \def\next{\errmessage{Use \string\Bbb\space only in math mode}}\fi\next}
\def\Bbb@#1{{\Bbb@@{#1}}}
\def\Bbb@@#1{\fam\msbfam#1}
\newfont{\twelvegoth}{eufm10 scaled\magstep1}
\newfont{\tengoth}{eufm10}
\newfont{\eightgoth}{eufm8}
\newfont{\sixgoth}{eufm6}
\newfam\gothfam
\textfont\gothfam=\twelvegoth
\scriptfont\gothfam=\eightgoth
\scriptscriptfont\gothfam=\sixgoth
\def\frak{\frak@}
\def\frak@#1{{\fam\gothfam{{#1}}}}
\def\frak@@#1{\fam\gothfam#1}
\catcode`@=12
%
%
%
\def\CC{{\Bbb C}}
\def\NN{{\Bbb N}}
\def\QQ{{\Bbb Q}}
\def\RR{{\Bbb R}}
\def\ZZ{{\Bbb Z}}
\def\cA{{\cal A}}          \def\cB{{\cal B}}          \def\cC{{\cal C}}
\def\cD{{\cal D}}          \def\cE{{\cal E}}          \def\cF{{\cal F}}
\def\cG{{\cal G}}          \def\cH{{\cal H}}          \def\cI{{\cal I}}
\def\cJ{{\cal J}}          \def\cK{{\cal K}}          \def\cL{{\cal L}} 
\def\cM{{\cal M}}          \def\cN{{\cal N}}          \def\cO{{\cal O}}
\def\cP{{\cal P}}          \def\cQ{{\cal Q}}          \def\cR{{\cal R}} 
\def\cS{{\cal S}}          \def\cT{{\cal T}}          \def\cU{{\cal U}}
\def\cV{{\cal V}}          \def\cW{{\cal W}}          \def\cX{{\cal X}}
\def\cY{{\cal Y}}          \def\cZ{{\cal Z}}
\def\qed{\hfill \rule{5pt}{5pt}}
\def\id{\mbox{id}}
\def\ggo{{\frak g}_{\bar 0}}
\def\uqggo{\cU_q({\frak g}_{\bar 0})}
\def\uqggp{\cU_q({\frak g}_+)}
\def\half{\frac{1}{2}}
\def\btf{\bigtriangleup}
\newtheorem{lemma}{Lemma}
\newtheorem{prop}{Proposition}
\newtheorem{theo}{Theorem}
\newtheorem{Defi}{Definition}

\vfill
\vfill
\begin{flushright}
math.QA/9807100
\end{flushright}
\begin{center}

{\LARGE {\bf {\sf 
Maps and twists relating $U(sl(2))$ and the nonstandard $U_{h}(sl(2))$: 
unified construction
}}} \\[0.8cm]

{\large B. Abdesselam$^{1}$\footnote{invite@elbahia.cerist.dz}, A.
Chakrabarti$^{2}$\footnote{chakra@orphee.polytechnique.fr}, 
R. Chakrabarti$^{3}$ and J. Segar$^{4}$}
\begin{center}
{\em 
$^{1}$Institute of Computer Systems, University of Mascara, Mamounia
Street,\\
2900 Mascara, Algeria. \\
$^{2}$Centre de Physique Th\'eorique\footnote{Laboratoire Propre 
du CNRS UPR A.0014}, Ecole Polytechnique, 91128 Palaiseau Cedex, France.\\
$^{3}$Department of Theoretical Physics, University of Madras, \\
Guindy Campus, 
Madras-600025, India. \\
$^{4}$Department of Physics, St. Peter's Engineering College, Avadi,
Madras-600 054, India. \\[2.8cm]
}
\end{center}

\end{center}

\smallskip

\smallskip 

\smallskip

\smallskip

\smallskip

\smallskip 

\begin{abstract}
A general construction is given for a class of invertible maps between the
classical $U(sl(2))$ 
and the Jordanian $U_{h}(sl(2))$ algebras. Different maps are directly useful 
in different contexts. Similarity trasformations connecting them, in so far as 
they can be explicitly constructed, enable us to translate results obtained in 
terms of one to the other cases.
Here the role of the maps is studied in the context of construction of twist
operators between the cocommutative and noncocommutative coproducts of 
the $U(sl(2))$ 
and $U_{h}(sl(2))$ algebras respectively. It is shown that a particular 
map called the 
`minimal twist map' implements the simplest twist given directly by the 
factorized
form of the ${\cal R}_{h}$-matrix of Ballesteros-Herranz. For other maps 
the twist has 
an additional factor obtainable in terms of the similarity 
transformation relating 
the map in question to the minimal one. The series in powers of $h$ for the
operator performing this 
transformation may be obtained up to some desired order, relatively easily. 
An explicit example is given for one particularly interesting case. 
Similarly the
classical and the Jordanian antipode maps may be interrelated by a similarity
transformation. For the `minimal twist map' the transforming operator is  
determined in a closed form.  

\end{abstract}

\vfill
\newpage 

\pagestyle{plain}

\sect{Introduction}
In the previous articles [1]-[3] we have introduced two different sets 
of invertible maps between the classical $sl(2)$ algebra
\begin{equation}
[J_{0} , J_{+} ] =\pm 2 J_{\pm}, \qquad  [ J_{+} , J_{-} ]  =J_{0}
\end{equation}
and the corresponding nonstandard Jordanian $h$-deformed $U_{h}(sl(2))$ 
algebra [4] \begin{equation}
[H,X] =\frac{2}{h} \sinh \, hX,  \quad
[H,Y] = -Y (\cosh \, hX ) - (\cosh\, hX) Y, \quad 
[X,Y] = H,
\end{equation}
where obviously in the $h \rightarrow 0$ limit, we have $(X,Y,H)\rightarrow
(J_{+},J_{-},J_{0})$.  
The map presented in [1], and subsequently  generalized in [2] to the
Jordanian deformation of the $so(4)$ algebra, maintains the 
diagonalization of the generator $H$:
\begin{equation}
H = J_{0}, \quad 
X = \frac{2}{h} arctanh (\frac{h}{2} J_{+}), \quad
Y  = (1 - (\frac{h J_{+}}{2})^2 )^{\frac{1}{2}} \,  J_{-}  \, (1 - (\frac{h
J_{+}}{2})^2)^{\frac{1}{2}}.
\end{equation}

The inverse of the map (1.3) is easily obtained [1]. The diagonalization of $H$ 
postulated above leads to attractive properties and easy construction of the 
irreducible representations of $U_{h}(sl(2))$ algebra and its finite 
dimensional  $R_{h}$-matrices [1]. The second map [3]
\begin{equation}
e^{\pm h X} = \pm h J_{+} + ( 1 + (hJ_{+})^2)^{\frac{1}{2}}, \,\,\,
H = (1 + (hJ_{+})^2)^{\frac{1}{2}} J_{0}, \,\,\,
Y = J_{-} + \frac{h^2}{4}  J_{+}  ( 1 - J_{0}^2) 
\end{equation}
plays a key role in the construction of the Jordanian $R_{h}$ matrices of
the $U_{h}(sl(2))$ algebra as contraction limits of the $R_{q}$ matrices of 
the standard
$q$-deformed $U_{q}(sl(2))$ algebra [5]. The standard $R_q$
matrices were transformed [3] by the operator
\begin{equation}
E_{q}(\eta {\cal  J}_{+} )  \otimes E_{q}(\eta {\cal J}_{+}),  
\end{equation}
where $({\cal J}_{\pm},{\cal J}_{0})$ are the generators of the 
$q$-deformed $U_{q}(sl(2))$ algebra [5] and  $\eta = \frac{h}{(q-1)}$. 
The deformed exponential in (1.5) reads 
\begin{equation}
E_{q}(x) = \sum_{n=0}^{\infty}\frac{x^n}{[n] !}, \qquad  [n] \equiv
\frac{q^n-q^{-n}}{q-q^{-1}}.
\end{equation}
The factor $\eta$  is {\it {singular}} in the $q \rightarrow 1$  limit, but the
singularities systematically cancel rendering the transformed $R_q$-matrices
{\it {finite}}. In the $q \rightarrow 1$ limit, the transformed $R_{q}$ 
matrices, on 
account of the map (1.4), reduce  to the finite dimensional represenations of
the Ballesteros-Herranz (B-H)
form [6] of the universal ${\cal {R}}_h$ matrix of the Jordanian algebra: 
\begin{equation}
{\cal R}_{h} = exp(-h X \otimes e^{hX} H) \,\, exp(e^{hX} \otimes hX).
\end{equation}
The universal ${\cal R}_{h}$ matrix in (1.7) satisfies $\sigma \circ 
{\cal R}_{h} =
({\cal R}_{h})^{-1}$; and consequently is triangular in nature. The flip 
operator
$\sigma$ mentioned above permutes vector spaces as follows: $\sigma(x 
\otimes y) = y \otimes x$.  
In [3] the case $\half \otimes j$ was treated fully. It was pointed out that the
treatment may be generalised not only to the general case $j_{1} \otimes 
j_{2}$ but
also to higher dimensional algebras. The last point was illustrated with the 
$U_{q}(sl(3))$.

Following Drinfield's arguments [7], the classical {\it {cocommutative}}
coalgebraic structure (denoted by the subscript {\it {c}}) of the $U(sl(2))$
algebra
\begin{equation}
\bigtriangleup_{c}(J_{\delta})=J_{\delta} \otimes 1 + 1 \otimes 
J_{\delta},\,\,\, 
S_{c}(J_{\delta}) = -J_{\delta},\,\,\, \epsilon_{c} (J_{\delta}) = 0
\quad (\delta=\pm,0)
\end{equation}
and the {\it {noncocommutative}} coalgebraic structure (denoted by the subscript
{\it {n}}) of the Jordanian $U_{h}(sl(2))$ algebra [4]
\begin{eqnarray}
&&\bigtriangleup_{n}(X)=X \otimes 1+1 \otimes X,\quad 
\bigtriangleup_n(Y)=Y\otimes  e^{h X}+e^{-h X} \otimes Y,\nonumber \\
&&\bigtriangleup_n(H)=H\otimes e^{h X}+e^{-h X}\otimes H,\nonumber \\ 
&&S_{n}(X)=-X, \quad S_{n}(Y)=- e^{h X} Y e^{- h X},  \quad S_{n}(H)=-e^{h X} H
e^{-h X}, \nonumber \\  
&&\epsilon_{n}(X)=\epsilon_{n}(Y)=\epsilon_{n}(H)=0,
\end{eqnarray}
may be related via suitable twist operators corresponding to various 
maps. Using the
maps (1.3) and (1.4) as examples, the relevant twist operators were 
considered in [3] as series expansions in $h$. 
 
The situation may be envisaged as follows. Different maps
may arise naturally in different contexts, and  may be particularly useful for 
different
purposes. As these maps relate the same pair of algebras, there exists an
equivalence relation between any two. The twist corresponding to 
different maps may
then also be related via the same relation. A unified treatment for a 
class of maps
may enable us to fully exploit the attractive properties of different 
maps in different
situations, and then obtain the results corresponding to others through 
equivalence relations, in so far as they can be obtained explicitly. This is 
attempted in the following sections.
In particular, we construct explicitly the `minimal twist map' for which 
the twist
operator $F$, relating the classical and the quantum Jordanian coproduct 
structures, corresponds directly to the factorized B-H form (1.7) of
${\cal R}_h$: namely,
\begin{equation}
F = exp(- e^{h X} H \otimes  hX).
\end{equation}
Then it is indicated how equivalence transformations can lead to the twists
corresponding to other maps. In the context of the above `minimal twist map',
the operator interrelating the classical and
the Jordanian antipode structures in (1.8) and (1.9) respectively may also be 
determined in a closed form.

In the following sections we often use the elements $T^{\pm 1} =
e^{\pm h X}$ of the $U_{h}(sl(2))$ algebra for calculational purpose.
In terms of these operators, the defining relations (1.2) and (1.9) of 
$U_{h}(sl(2))$ algebra read
\begin{equation}
[H , T^{\pm 1 }] = T^{\pm 2} - 1,\,\,\, [H, Y] = -\half ( Y(T+T^{-1}) +
(T+T^{-1} )Y ),\,\,\, [T^{\pm 1} , Y ] = \pm \frac{h}{2} (H T^{\pm 1} +
T^{\pm 1} H)
\end{equation}
and 
\begin{eqnarray}
&&\bigtriangleup_{n} (T^{\pm 1}) = T^{\pm 1} \otimes  T^{\pm 1}, \,\,
\bigtriangleup_n(Y) = Y \otimes T + T^{-1} \otimes Y, \,\,
\bigtriangleup_{n}(H) = H \otimes T + T^{-1} \otimes H, \nonumber \\
&& S_{n}(T^{\pm 1}) = T^{\mp 1} , \,\,\, S_{n}(Y) = -TYT^{-1}, \,\,\, S_{n}(H) 
= -THT^{-1},\nonumber \\
&&\epsilon_{n}(T^{\pm 1} ) = 1, \,\,\,\epsilon_{n}(Y) = \epsilon_{n}(H) =
0.
\end{eqnarray}
In the constructions to follow all functions of $J_{+}$ or $X$ may be 
interpreted as a
finite power series on spaces of irreducible representations. This provides well
defined operators.
\sect{A class of maps}
To construct a class of maps important for our purpose, we take an ansatz 
 \begin{equation}
J_{+} = f_{1}, \,\,\,
J_{0} = f_{2} H, \,\,\, 
J_{-} = f_{3} Y + u + v H + w H^2,
\end{equation}
where we introduce $(f_1, f_2, f_3; u,v, w)$ as functions of $T$ only. 
We start with
(2.1) rather than with (2.14) to follow since, in the context of the 
former, the map
crucial for our determination of the twist operator assumes a simple 
form. An additive function $f(T)$ in the expression 
for $J_{0}$ may be absorbed by a similarity transformation. 
To ensure correct classical limits, the introduced functions are required to 
satisfy the limiting properties  
\begin{equation}
(f_1,f_2,f_3;u,v,w) \rightarrow (X,1,1;0,0,0)
\end{equation} 
as $h \rightarrow 0$. We will, moreover, be interested in invertible maps.
For any function $f(T)$, we denote 
$f^{\prime} = \frac{d}{dT}f$. Using the identities 
\begin{eqnarray}
&&[H,f] = (T^2-1) f^{\prime}, \nonumber \\
&&[f,Y] = \frac{h}{2}\left( (T f^{\prime} )H + H(T f^{\prime}) 
\right)=h(T f^{\prime} )H +\frac{h}{2}(T^2-1) (Tf^{\prime})^{\prime}
\end{eqnarray} 
and the algebraic constraints (1.1) and (1.11) systematically, we, {\it 
{for a given $f_1$}}, obtain  a set of {\it {seven}} coupled nonlinear equations
(one being of second order in
$f_2$) for the {\it {five}} unknown functions $(f_2 , f_3;u,v,w)$: 
\begin{eqnarray}
&&(T^2-1)f_2f_{1}^{\prime}=2f_1,\,\,\,hTf_3f_1^{\prime}-2w(T^2-1)f_1^{\prime}
=f_2, \nonumber \\
&&\frac{h}{2}f_3(Tf_1^{\prime})^{\prime}-vf_1^{\prime}-
w((T^2-1)f_1^{\prime})^{\prime}=0,\,\,\,
f_2 \left( (T^2-1)f_3^{\prime}-(T+T^{-1})f_3
\right)=-2f_3,\nonumber \\
&&(T^2-1)f_2 \left( \frac{h}{4}(1+T^{-2})f_3+u^{\prime}\right)=-2u,\,\,\,
hTf_3f_2^{\prime}+(T^2-1)(f_2w^{\prime}-2wf_2^{\prime})=-2w, \nonumber \\
&&(T^2-1)\left( \frac{h}{2}f_3(f_2+T(Tf_2^{\prime})^{\prime})T^{-1}
-(vf_2^{\prime}-f_2v^{\prime})-w((T^2-1)f_2^{\prime})^{\prime} \right) =-2v. 
\end{eqnarray} 
These equations may  then be 
solved consistently and the limiting behaviour (2.2) as $h \rightarrow 
0$ may be 
implemented through proper choice of integration  constants. This leads 
unambiugously to the following solution
\begin{eqnarray}
&&f_2 = \frac{2}{(T^2-1)} \frac{f_1}{f_{1}^{\prime}}, \,\,\,\,
f_3 = \frac{1}{2h} \frac{(T- T^{-1})}{f_1}, \,\,\,\,
u = - \frac{(T - T^{-1})^2}{16 f_1}, \nonumber \\
&&v = -\half \left( \frac{f_2^\prime}{f_1^\prime} - \frac{f_2}{f_1} 
+ \frac{T + T^{-1}}{2f_1} \right), \,\,\,\,
w = \frac{1}{4f_1}(1 - f_2^2).
\end{eqnarray}
In the solution (2.5) the unknown functions ($f_2,f_3; u,v,w$) may be fully
expressed in terms of the function $f_1$ and its derivatives, assumed to be
known. Thus a suitably chosen $f_1$ and its
properties completely determine the map. The simplest choice
\begin{equation}
f_1 = \pm \frac{1}{h} (T^{\pm 1} - 1)
\end{equation}
do not seem to lead to particularly interesting properties. This is true for our
present main concern, namely construction of twists. But (2.6) probably deserves
further study in other contexts.
The inverse of the map (1.3) corresponds to the choice $f_1 = \frac{2}{h}$
$\left( \frac{T-1}{ T+1}\right)$. The unknown functions now, through the
solution (2.5), assume the form  
\begin{eqnarray}
&&f_2 =1,\,\,\, f_3=\left( \frac{T^{\half} +T^{-\half}}{2} \right)^2,\,\,\,
u=-\frac{h}{8}\left( \frac{T^\half+T^{-\half}}{2} \right)^2
(T-T^{-1}),\nonumber \\
&&v=-\frac{h}{8}(T-T^{-1}),\,\,\,\,
w=0.
\end{eqnarray}
The choice $f_1 = \frac{1}{2h}(T - T^{-1})$ reproduces the inverse of the
map (1.4). The solution (2.5) now restricts the unknown functions as 
\begin{eqnarray}
&&f_2=\frac{2}{T+T^{-1}},\,\,\, f_3=1,\,\,\, u = -\frac{h}{8}(T -
T^{-1}),\,\,\, \nonumber \\
&&v = -\frac{h}{2} \left( \frac{T-T^{-1}}{T+T^{-1}} \right)^3 
,\,\,\, w = \frac{h}{2} \frac{( T-T^{-1})}{(T + T^{-1})^2}.
\end{eqnarray}
The special interest of (2.7) and (2.8) have already been indicated. 
Other choices of
$f_1$ may prove interesting. For a purpose of this article, namely 
finding the map corresponding to the form (1.10) of the twist
operator, the pertinent choice is 
\begin{equation}
f_1 = \frac{1}{2h}(1-T^{-2}). \nonumber \\
\end{equation}
The solution (2.5) now yields 
\begin{equation}
f_2=f_3=T,\,\,\, u = -\frac{h}{8}(T^2-1),\,\,\, v = 
-\frac{h}{2}(T^2-1)T,\,\,\,w= -\frac{h}{2}T^2. 
\end{equation}
This is the `minimal twist map' resulting in the twist operator (1.10) 
as will be shown in a subsequent section.

Even for invertible maps it is interesting to construct the general 
solution starting
from the other end. So let us consider directly the class of maps given by 
\begin{equation}
T=g_1,\,\,\, 
H=g_2 J_{0},\,\,\, 
Y=g_3 J_{-} + a + b J_{0} + c J_{0}^2, 
\end{equation}
where $(g_1,g_2,g_3; a,b,c)$ are functions of $J_{+}$ only. The expansion 
\begin{equation}
g_1 = 1 + h J_{+} + O(h^2). 
\end{equation}
conforms to the correct limiting properties as $h \rightarrow 0$; and, 
in the same limit, the other functions behave as
\begin{equation}
(g_2,g_3;a,b,c) \rightarrow (1,1;0,0,0).
\end{equation}
For any function $g(J_{+})$, we denote its derivative as $g^{\prime} = \frac{d
g}{dJ_{+}}$. The following identities  
\begin{equation}
[J_{0} , g] = 2 J_{+} g^\prime, \qquad
[g, J_{-} ] = \half (J_{0} g^\prime + g^\prime J_{0} )
 = J_{+} g^{\prime \prime}  + g^\prime J_{0}
\end{equation}
and a systematic use of the algebraic properties (1.1) and (1.11) yield, 
as before, 
a set of {\it {seven}} coupled nonlinear equations. These equations, 
after taking into 
account the limiting properties (2.11) and (2.12) for determination of 
the constants 
of integrations, may be solved unambiguously.  The discussion following 
(2.3) is also relevant here.  The unknown functions read 
\begin{eqnarray}
&&g_{2} = \frac{\left(g_1^2 - 1 \right)}{2 J_{+} g_1^\prime},\,\,\,
g_3=2h \frac{J_{+} g_1}{\left( g_1 ^2 - 1 \right)},\,\,\,
a= \frac{h}{8}(g_1-g_1^{-1}), \nonumber \\
&&b = - \left( 2c + \frac{h}{2} (\frac{g_1}{g_1^\prime})^\prime g_2 \right)
,\,\,\,c = \frac{h g_1(1-g_2^2)}{2(g_1^2-1)}.
\end{eqnarray} 
Thus again, as is to be expected, the choice of $g_1$ and the limiting 
properties
as $h \rightarrow 0$ completely determine the map. The three particular 
invertible solutions introduced before well illustrate the situation.
\sect{Equivalence of the maps}
As emphasized before, different maps play useful roles in specific 
contexts. These
maps are, however, equivalent to one another through appropriate similarity 
transformations. Understanding the equivalence properties precisely is 
of importance 
as they allow us to carry over the results easily derived using one map 
in the context 
of another one. This will be illustrated in the construction of twists.
The class of maps  (2.11) discussed above, may be related through similarity 
transformations by operator of the form $exp( \lambda (J_{+}) J_{0})$. 
Alternately,
the class of inverse maps (2.1) may be interrelated in a parallel way by 
the operator 
$exp(\mu(T)H)$.  In the classical $h \rightarrow 0 $ limit, the 
generating functions 
$\lambda(J_{+})$ and $\mu(T)$ introduced above behave as $(\lambda, \mu 
)\rightarrow
0$. It is usually difficult to obtain the exact, closed form expressions
of the generating functions $\lambda$ and $\mu$. But
the constant coefficients in the series for, say,
\begin{equation}
\lambda(J_{+}) = c_1(h J_{+} ) + c_2 (h J_{+})^2 + \cdot \cdot \cdot + 
c_{n} (hJ_{+})^n + \cdot \cdot \cdot,
\end{equation}
may be obtained iteratively in a systematic fashion. We will demonstrate 
this in the
following. A similar series  in $(h X)$ for the function $\mu(T)$ may also be
constructed. 

Let us consider two maps starting respectively with
\begin{equation}
T = g_{1}(J_{+})\,\,\, and
\,\,\, {\hat T} = {\hat g_{1}}(J_{+}).
\end{equation}
To avoid confusion we distinguish between the two sets of generators 
$(T,H,Y)$ 
and (${\hat T}$, ${\hat H}$, ${\hat Y}$), both satisfying (1.11) 
corresponding to the
two maps. Our task now is to construct the generating function 
$\lambda(J_{+})$ that transforms one map to another as 
\begin{equation}
e^{- \lambda (J_{+}) J_{0} } g_1(J_{+}) e^{\lambda (J_{+}) J_{0}} = {\hat
g}_{1}(J_{+}).
\end{equation}
It is important to realize that {\it {it is sufficient to ensure}} 
(3.3). The required
equivalence of the pairs $(H,{\hat H})$ and $(Y, {\hat Y})$ follow  
once that of $(T, {\hat T})$ is assured. This is an evident consequence of the 
following facts: \\
(1) A well defined similarity transformation conserves the algebra.  \\  
(2) The transforming operator in (3.3) preserves the form of the ansatz 
(2.11) altering only the functions of $J_{+}$. \\
(3) For our maps and limiting $(h \rightarrow 0 )$ constraints, a choice of the
function $g_1 ({\hat g}_{1})$  determines the remaining functions 
unambiguously. \\ The equivalence relation (3.3) may be recast in the form
\begin{equation}
g_{1} \left( e^{-\lambda(J_{+}) J_{0}} \, J_{+} \, e^{\lambda(J_{+}) 
J_{0}} \right) = {\hat g}_{1}(J_{+}).
\end{equation}
The invertibility of our maps, in conjunction with (3.4), now yield
\begin{equation}
e^{-\lambda(J_{+}) J_{0}} \, J_{+} \,  e^{\lambda(J_{+}) J_{0}} = f_{1} 
\left( {\hat g}_{1}(J_{+}) \right).
\end{equation}
The transforming relation (3.5) may be systematically used to 
generate the series (3.1), 
where the coefficients $\{c_{i}, |i=1,2,... \}$ may be evaluated 
iteratively. 
This is best illustrated through a relatively simple but particularly 
interesting 
example. We consider the maps (1.4) and (2.9) along with their inverses. 
The relation (3.5) now implies 
\begin{equation}
e^{- \lambda(J_{+}) J_{0} } (h J_{+} ) e^{\lambda(J_{+}) J_{0} } = (h J_{+}) 
\left( -(hJ_{+}) +  \left( 1 + (hJ_{+})^2 \right)^{\half} \right).
\end{equation}
Using standard expansion scheme in (3.6), we now obtain the generating
function $\lambda(J_{+})$ up to $O(h^5)$:
\begin{equation}
\lambda(J_{+}) = \half(hJ_{+}) +
\frac{1}{4}(hJ_{+})^2+\frac{1}{8}(hJ_{+})^3+\frac{1}{24}(hJ_{+})^4
-\frac{1}{96}(hJ_{+})^5+\cdot\cdot\cdot
\end{equation}
The higher order terms may be computed similarly. 
Using our maps we may express the transforming operator 
$exp{\left(  \lambda (J_{+})J_{0} \right)}$ as $exp{\left( \mu (T) 
H\right)} $ and vice versa.
An application of such equivalence will be indicated in the following 
section. Unlike the general solutions
(2.5) and (2.15) for the class of maps considered here, we are, however, 
unable to 
present an exact general solution for (3.3) in an explicit form.
\sect{Twists: the role of maps}
The classical algebra (1.1) with the cocommutative coalgebra structure (1.8)
may be assumed to have a trivial universal ${\cal R}$ matrix: 
\begin{equation}
{\cal R}_{c} = 1 \otimes 1.
\end{equation}
On the other hand, the factorized B-H form [6] of the universal ${\cal 
R}_{h}$ matrix (1.7) of the $U_{h}(sl(2))$ algebra may be recast as 
\begin{equation}
{\cal R}_{h} = (\sigma \circ V) V^{-1}, 
\end{equation}
where $V^{\pm 1} = exp(\mp e^{h X} \otimes h X) = exp(\mp T H \otimes h X)$. 
The flip operator $\sigma$ has been introduced earlier following equation (1.7).
Following [7] it may be observed that the triangular universal ${\cal 
R}_{h}$ matrix 
(1.7) may, by a suitable twist  $F \in U_{h}(sl(2))^{\otimes 2}$, be 
brought into
the classical form (4.1). Such twists, relating the quantum comultiplication 
(1.9) with classical cocommutative coproduct (1.8), may be called factorizing
twist as they factorize [7] the universal ${\cal R}_{h}$ matrix as 
\begin{equation}
{\cal R}_{h} = (\sigma \circ F ) \, {\cal R}_{c} \, F^{-1} = (\sigma 
\circ F) F^{-1}. \end{equation}
The second equality in (4.3) follows from (4.1). The consistency of the 
results (4.2) and (4.3) requires 
\begin{equation}
F = V F_S,
\end{equation}
where $F_S$ is symmetric under permutation: $\sigma \circ F_{S} = F_{S}$.
In the structure (4.4) of the twist operator $F$ the factor $V$ is fixed;
whereas different choices of $F_{S}$ imply different maps between the relevant
algebras.
For an invertible map 
\begin{equation}
m: (T^{\pm 1}, H,Y) \rightarrow (J_{\pm},J_{0}),\qquad 
m^{-1}: (J_{\pm},J_{0}) \rightarrow (T^{\pm 1},H,Y)
\end{equation}
the coproducts given by (1.8)
and (1.9) respectively may be related via the twist operator $F$ as 
\begin{equation}
F^{-1}  (  \btf_{n} (\Phi ) )  F = {\tilde {\bigtriangleup}}_{c} (\Phi),  \quad 
\forall  \Phi \in U_{h}(sl(2)),
\end{equation}
where, ${\tilde {\bigtriangleup}}_{c} = (m^{-1} \otimes m^{-1} ) 
\circ \btf_{c} \circ m$. As a consistency check of our factorization 
scheme (4.4) 
for the twist operator $F$, we utilize (4.6) to obtain the identity
\begin{eqnarray}
F_{S} ( {\tilde {\bigtriangleup}}_{c} (\Phi) ) F_{S}^{-1} &=& V^{-1} 
( \bigtriangleup_{n} (\Phi) ) V \nonumber \\
& =& \half  V^{-1} \left( \btf_{n} (\Phi ) + {\cal R}_{h}^{-1}(\sigma 
\circ \btf_n (\Phi)){\cal R}_{h} \right)V  \nonumber \\
&=& \half  \left(  V^{-1} ( \btf_{n} (\Phi) ) V + \sigma \circ \left( 
V^{-1} ( \btf_{n} (\Phi) ) V \right) \right),
\end{eqnarray}
where we have used the standard properties of the universal ${\cal 
R}_{h}$ matrix and
its factorized form (4.2). As the operator ${\tilde \btf}_{c}(\Phi)$ and 
$F_{S}$ on the
lhs of (4.7) are both symmetric under the permutation, the identity 
explicitly shows the
consistency of our factorization (4.4) of the twist operator $F$ with 
the defining relation (4.6). 
At this stage it is natural to ask the following question:
Is it possible to  consistently set $F_{S} = 1 \otimes 1 ?$ And if so, 
what does 
it imply? An affirmative answer to the above question, in view of (4.4), 
leads to
the following elegant and useful construction of the twist operator: 
\begin{equation}
F = V \left( = exp(-e^{h X} H \otimes h X) \right). 
\end{equation}
We show that (4.8)  can indeed be realized and that it implies the
implementation of a specific map.
This map will be constructed explicitly, step by step, by exploring the 
action of $V$
on suitable operators. The following results are crucial. Setting $V^{-1}
=exp(-\Lambda)$, we assume the construction (4.8) to be vaild and thereby  
evaluate the lhs of (4.6) using the standard Campbell-Hausdorff series 
expansion: \begin{equation}
V^{-1} ( \btf_{n} (\Phi) ) V = \btf_{n} (\Phi) + \left[ \Lambda , \btf_{n} 
(\Phi )\right] + \half
\left[ \Lambda , \left[ \Lambda , \btf_n (\Phi) \right] \right] + \cdot 
\cdot \cdot .
\end{equation}
Our task is now to make appropriate choices for $\Phi$ to realize 
the map, thus  validating the construction (4.8) of the twist operator $F$. 
For the choice 
$\Phi = (1 - T^{-2}) $,  the coproduct rules (1.12) yield
\begin{equation}
\btf_{n} (\Phi) = 1 \otimes 1 - T^{-2} \otimes T^{-2}. 
\end{equation}
Using (4.8)-(4.10) we obtain a remarkable result
\begin{equation}
V^{-1} \left( \btf_{n} \left( 1 - T^{-2} \right) \right) V = \left(1 - 
T^{-2} \right) \otimes 1 + 1 \otimes \left(1 - T^{-2} \right).
\end{equation}
Similarly the following result  
\begin{equation}
\btf_{n}(T H ) = T H \otimes T^2 + 1 \otimes T H
\end{equation}
and the similarity transformation (4.9) yield 
\begin{equation}
V^{-1} \btf_{n} \left(TH \right) V  = T H \otimes 1 + 1 \otimes T H.
\end{equation}
The `classical' aspect of the r.h.s. of (4.11) and (4.13) 
is not an accident. The map (2.9) and (2.10) may be written 
after a regrouping of the terms for $J_{-}$ as
\begin{equation}
J_{+}= \frac{1}{2h} \left( 1 - T^{-2} \right),\,\,
J_{0} = T H,\,\,\
J_{-} = T Y - \frac{h}{2} (T H)^2 - \frac{h}{8} ( T^2 - 1). 
\end{equation}
Indeed it may be demonstrated by direct computation that the expression for
$J_{-}$ satisfies the requirement (4.6):
\begin{equation}
V^{-1} \left( \btf_{n} (J_{-}) \right) V = J_{-} \otimes 1 + 1 \otimes J_{-} 
(= {\tilde \btf}_{c} (J_{-}) )
\end{equation}
The equation (20) of [6] is quite close to this result. But the authors, 
their context being different, did not construct the map completely.
The above results (4.11)-(4.15) may be combined as 
\begin{equation}
V^{-1} (\btf_{n} (J_{\delta}) ) V = J_{\delta} \otimes 1 + 1 \otimes
J_{\delta} (={\tilde \btf}_{c} (J_{\delta}) ), \quad (\delta=\pm,0).
\end{equation}
In (4.15) and (4.16) the classical generator $J_{\delta}$ stand for the 
expressions in the rhs of inverse map (4.14).
Comparing (4.16) with (4.7), we now obtain 
\begin{equation}
F_{S} ( {\tilde \btf}_{c} (J_{\delta}) ) F_{S}^{-1} = {\tilde \btf}_{c} 
(J_{\delta}).
\end{equation}
This verifies our previously announced structure of the twist operator
corresponding to the inverse map (4.14):
\begin{equation}
F_{S} = 1 \otimes 1 \quad  \Rightarrow  \quad F= V.
\end{equation}
From this point of view we call (4.14) the `minimal  twist map'. The twist
operator $V$ in (4.8) satisfies the cocycle condition:
\begin{equation}
\left( (\btf_n \otimes 1) V \right ) (V \otimes 1) = \left( (1
\otimes \btf_n ) V \right ) (1 \otimes V).
\end{equation} 

Starting with the `minimal twist map', $F_{S}$ corresponding to another 
one may be expressed in terms of the similarity transformation relating them.
Let $J_{\delta} \,\,\, (\delta=\pm,0)$ correspond to the `minimal twist map'
(4.14) and ${\hat J}_{\delta}$ to another
map related to the former by a transformation $U$ of the class discussed 
before in the context of  (3.2)-(3.7), such that 
\begin{equation}
J_{\delta} = U {\hat J}_{\delta} U^{-1}. 
\end{equation}
The equivalence property (4.20) and the following identity, 
\begin{equation}
F_{S} ({\tilde \btf}_{c} ({\hat J}_{\delta}) ) F_{S}^{-1} = V^{-1} 
\left( \btf_{n} ({\hat J}_{\delta})\right) V
\end{equation}
obtained {\` a} la (4.7) now yield 
\begin{equation}
F_{S} \left( {\tilde \btf}_{c} U \right)^{-1} \left( {\tilde \btf}_{c} 
J_{\delta} \right) \left( {\tilde \btf}_{c} U \right) F_{S}^{-1} = 
\left(V^{-1} \left( \btf_{n} U^{-1} \right) V \right) \left( V^{-1} 
\left(\btf_{n} J_{\delta} \right) V \right) \left( V^{-1} \left(
\btf_{n} U \right) V \right).
\end{equation}
The property (4.16) of the `minimal twist map'  $J_{\delta}$ now imply 
\begin{equation}
F_{S} = \left( V^{-1} \left(\btf_{n} (U) \right) V\right)^{-1} 
\left({\tilde \btf}_{c} (U) \right), 
\end{equation}
where not only the factor $\left({\tilde \btf}_{c} (U)\right)$ on the rhs but
first factor, as follows from (4.7), is also symmetric under permutation. 
This assures consistency with the 
postulated symmetry  of $F_{S}$. Finally, using (4.4) and (4.23), we obtain the
general structure of the twist operator for arbitrary maps of the class (2.1):
\begin{equation}
F = \left( \btf_{n} (U) \right)^{-1} V \left( {\tilde \btf}_{c} (U) \right).
\end{equation}
In (4.20)-(4.24) the equivalence operator $U$ is of the form $ U = 
exp(-\mu(T)H)$ 
discussed earlier in section 3. The twist operators corresponding to 
any two maps may now be related in an evident way.

An example of $U$ is considered in (3.6) and its series development is 
indicated 
in (3.7). If instead of the map (1.4) we consider (1.3), then (3.6) is to be 
replaced by
\begin{equation}
e^{-\lambda(J_{+}) J_{0}} (h J_{+} ) e^{\lambda(J_{+}) J_{0}} = 
\frac{h J_{+} }{\left( 1 + \frac{h J_{+}}{2} \right)^2}
\end{equation}
and the coeffiecents $c_{n}$ in (3.1) must be evaluated accordingly. For the 
map (1.3), starting with the coproduct $\btf_{n}$ in (1.9), it is possible to 
write down the induced coproducts for
$J_{\delta}$. An interesting application of these was made [8]-[10] 
in computing the
C.G. coeffiecents of $U_{h}(sl(2))$. From our point of view, as per (4.6), 
such induced coproducts
$\btf_{n} (J_{\delta})$ may be finally expressed in terms of the relevant 
classical coproduct and the twist operator (4.24). The cocycle condition for an 
arbitrary twist operator $F$ in (4.24) may be verified by using an appropriate 
series expansion for $U$. 

The antipode structure of Jordanian $U_{h}(sl(2))$ algebra (1.12) and that
of the classical $U_{h}(sl(2))$ algebra (1.8) may also be related by a 
similarity transformation as follows:
\begin{eqnarray}
G^{-1} S_{h} (\Phi) G = m^{-1} \circ S_{c} \circ m (\Phi)\quad \Leftrightarrow 
\quad m
\circ S_{h} \circ m^{-1} ({\tilde \Phi}) = {\tilde G} S_{c} ({\tilde \Phi})
{\tilde G}^{-1}, \nonumber \\ 
\forall \Phi \in U_{h}(sl(2)),\,\, {\tilde \Phi} \in U(sl(2)),
\end{eqnarray}
where $G \in U_{h}(sl(2))$ and ${\tilde G} = m \circ G \in U(sl(2))$. 
For the `minimal twist map' (4.14), the operator $G$ and ${\tilde G}$ 
performing the
similarity transformation (4.26) may be determined in a closed form:
\begin{equation}
G = \exp \left( \half T H (1 - T^{-2}) \right) \quad \Leftrightarrow \quad 
{\tilde G} = \exp (h J_{0} J_{+}).
\end{equation}
For the other maps discussed in section 2, these transforming operators may 
be determined as power serieses in the deformation parameter $h$. We will
not discuss this here. The counit maps in (1.8) and (1.9) readily 
correspond to each other.   
\sect{Remarks}
Contrasting an earlier approach [3]  of starting with the `classical' 
$r$-matrix and 
then trying to incorporate higher order terms in the construction  of 
a series expansion of the twist 
operator corresponding to a specific map, here we start from the very 
simple fact that
for a classical algebra the ${\cal R}$-matrix is identity. Then we note that 
the B-H form of ${\cal R}_{h}$ 
obtained in [6] is already expressed as the twisted form of the identity 
and the 
margin for manoeuvre left, as evidenced from (4.4) to (4.7), is a factor 
$F_{S}$ 
symmetric under permutation. We then show that constraining this factor 
to be unity 
implies a particular invertible map relating the generators of $U(sl(2))$ and 
$U_{h}(sl(2))$ algebras. We call this the `minimal twist  map'. But 
other maps are of 
interest in other contexts. So we give a unified  construction of a class of 
invertible maps, their equivalence relations and the constructions of $F_{S}$ 
in terms of the latter ones. The simplest possible twist is undoubtedly given 
by (4.8), which corresponds to the `minimal twist map' (4.14).
For the `non-minimal' cases it is much simpler to compute $U$, the similarity
transformation to the minimal case and then $F_{S}$ in terms of $U$ as in
(4.23). An example of series development for $U$ is given in (3.7). 
Combining the
above results in (4.24) we obtain, for an arbitrary map of the class 
discussed in
section 3, the twist operator interrelating the classical and the 
noncocommutative Jordanian coproducts.  The classical and
the Jordanian antipode maps may be related through a similarity transformation. 
The transforming 
operator corresponding to `minimal twist map' may be obtained in a 
closed form. For an arbitrary map, the pertinent transforming 
operator may be obtained as a series in $h$.

Other, quite different looking, constructions for ${\cal R}_{h}$ [11] 
and the twist
operator [4] may be found. We will not attempt to trace their 
possible relations to the formalism presented here. Finally we just briefly 
mention that after the first twist leading from the cocommutative to the 
noncocommutative 
Hopf structure, it is possible to envisage a second twist leading 
to a quasi Hopf generalization 
for the nonstandard case. This will be studied in a following paper.

\bibliographystyle{amsplain}

\end{document}